\documentclass[12pt]{article}
\def\hind{\hangindent=2pc\hangafter=1}
\newfont{\smcaps}{cmcsc10 scaled\magstep1}

\newcommand{\Cov}{{\rm ~Cov\,}}
\newcommand{\E}{{\rm ~E\,}}
\newcommand{\LCOV}{{\rm ~LCOV\,}}
 at 11truept
\usepackage{epsfig}
\raggedright
\begin{document}
\baselineskip=22pt

\title{COMPUTER ALGEBRA DERIVATION OF THE BIAS OF BURG ESTIMATORS}
\date{September 14, 2004}
\author{Y. ZHANG and A.I. MCLEOD\\
Department of Statistical and Actuarial Sciences,\\
The University of Western Ontario, London, Ontario N6A 5B7 \\
Canada
}

%\author{Author}
\maketitle
%\newpage
\hrule
\bigskip
{\bf Abstract.\/}
A symbolic method is discussed which can be used to obtain the asymptotic
bias and variance to order $O(1/n)$ for estimators in stationary time series.
Using this method the bias to $O(1/n)$ of the Burg estimator in AR(1) and
AR(2) models is shown to be equal to that of the least squares
estimators in both the known and unknown mean cases. Previous
researchers have only been able to obtain simulation results for
this bias because this problem is too intractable without using
computer algebra.

\bigskip

\noindent \underline{Preprint}\\
Ying Zhang and A. Ian McLeod (2006), Computer Algebra Derivation of the Bias of Burg Estimators, {\it Journal of Time Series Analysis} 27, 157-165

\newpage
\begin{center}
1. INTRODUCTION AND SUMMARY
\end{center}
\medskip
Tj\o stheim and Paulsen (1983, Correction 1984) showed that the
Yule-Walker estimates had very large mean-square errors in the
AR(2) case when the parameters were near the admissible boundary
and that this inflated mean square error was due to bias. This
result was demonstrated by Tj\o stheim and Paulsen (1983) in
simulation experiments as well as by deriving the theoretical bias
to order $O(1/n)$. It was also mentioned by Tj\o stheim and
Paulsen (1983, p.397, \S 5) that the bias results from simulation
experiments for the Burg estimates were similar to those obtained
for least squares estimates but that they had not been able to
obtain the theoretical bias term. Using computer algebra we are
now able to compute the bias for the Burg estimator and to show
that to order $O(1/n)$ it is equal to the bias for the least
squares estimator in the case of AR(1) and AR(2) models in both
the known and unknown mean cases.

As pointed out by Lysne and Tj\o stheim (1987), the Burg
estimators have an important advantage over the least squares
estimates for autoregressive spectral estimation since Burg
estimates always lie in the admissible parameter space whereas the
least squares estimates do not.  Burg estimators are now
frequently used in autoregressive spectral estimation (Percival
and Walden, 1993, \S 9.5) since they provide better resolution of
sharp spectral peaks. As the Yule-Walker estimators, the Burg
estimators may be efficiently computed using the Durbin-Levinson
recursion.  Our result provides further justification for the
recommendation to use the Burg estimator for autoregressive
spectral density estimation as well as for other autoregressive
estimation applications.

It has been shown that symbolic algebra could greatly simplify
derivations of asymptotic expansions in the IID case (Andrews and
Stafford, 1993). Symbolic computation is a powerful tool for
handling complicated algebraic problems that arise with expansions
of various types of statistics and estimators (Andrews and
Stafford, 2000) as well as for exact maximum likelihood
computation (Currie, 1995; Rose and Smith, 2000). Cook and
Broemeling (1995) show how symbolic computation can be used in
Bayesian time series analysis. Smith and Field (2001) described a
symbolic operator which calculates the joint cumulants of the
linear combinations of products of discrete Fourier transforms. A
symbolic computation approach to mathematical statistics is
discussed by Rose and Smith (2002). In the following sections, we
develop a symbolic computation method that can be used to solve a
wide variety of time problems involving linear time series
estimators for stationary time series. Our symbolic approach is
used to derive the theoretical bias to $O(1/n)$ of the Burg
estimator for the AR(2) model.

\bigskip
\begin{center}
2. ASYMPTOTIC EXPECTATIONS AND COVARIANCES
\end{center}

\medskip
Consider $n$ consecutive observations from
a stationary time series, $z_t, t=1,...,n$, with mean $\mu=\E(z_t)$ and autocovariance
function $\gamma_k = \Cov(z_t,z_{t-k})$.
If the mean is known, it may, without loss of generality be taken to be
zero. Then one of the unbiased estimators of autocovariance
$\gamma (m-k)$ may be written as
\begin{equation}
{\cal S}_{m,k,i} = {1\over n+1-i} \sum\limits_{t=i}^n z_{t-m}
z_{t-k},
\end{equation}
where $m$, $k$ and $i$ are non-negative integers with $\max (m,k) <
i \le n$. If the mean is unknown, a biased estimator of $\gamma
(m-k)$ may be written as
\begin{equation}
{\overline {\cal S}}_{m,k,i} = {1\over n} \sum\limits_{t=i}^n
(z_{t-m}-\overline {z}_n) (z_{t-k}-\overline {z}_n),
\end{equation}
where ${\overline z}_n$ is the sample mean.

{\smcaps Theorem 1.\/} Let the time series ${z_t}$ be the
two-sided moving average,
\begin{equation}
z_t=\sum_{ j=-\infty }^{\infty} \alpha_j e_{t-j},
\end{equation}
where the sequence $\{\alpha_j\}$ is absolutely summable and the
$e_t$ are independent $N(0,{\sigma}^2)$ random variables.  Then
for $i\le j$,
\begin{equation}
\lim_{n \to \infty} n \Cov( {\cal S}_{m,k,i}, {\cal S}_{f,g,j} ) =
\sum_{ h=-\infty }^{\infty} T_h,
\end{equation}
where $T_h = \gamma(g-k+i-j+h)
\gamma(f-m+i-j+h)+\gamma(f-k+i-j+h)\gamma(g-m+i-j+h)$.

{\smcaps Theorem 2.\/} Let a time series $\{z_t\}$ satisfy
the assumptions of Theorem 1.  Then
\begin{equation}
\lim_{n \to \infty} n \E( \overline {\cal S}_{m,k,i}-\gamma (m-k))
= -|i-1|\gamma (m-k)-\sum_{ h=-\infty }^{\infty} \gamma(h)
\end{equation}
and
\begin{equation}
\lim_{n \to \infty} n \Cov( \overline {\cal S}_{m,k,i}, \overline
{\cal S}_{f,g,j} ) = \sum_{ h=-\infty }^{\infty} T_h,
\end{equation}
 where $T_h = \gamma(g-k+i-j+h)
\gamma(f-m+i-j+h)+\gamma(f-k+i-j+h)\gamma(g-m+i-j+h)$.

These two theorems may be considered as the extensions of Theorem
6.2.1 and Theorem 6.2.2 of Fuller (1996). Letting $p=m-k$ and
$q=f-g$, the left side of (4) or (6) can be simplified,
\begin{equation}
\sum_{ h=-\infty }^{\infty} T_h =
\sum_{h=-\infty}^{\infty}\gamma(h)
\gamma(h-p+q)+\gamma(h+q)\gamma(h-p).
\end{equation}
There is a wide variety of estimators which can be written as a
function of the autocovariance estimators, ${\cal S}_{m,k,i}$ or
$\overline {\cal S}_{m,k,i}$, such as, autocorrelation estimator,
least squares estimator, Yule-Walker estimator, Burg estimator,
etc. The asymptotic bias and variance may be obtained by the
Taylor expansion.  Unfortunately, in the most cases, those
expansions include a large number of expectations and covariances
of the autocovariance estimators.  It is too intractable manually.
Theorems 1 and 2 provide the basis for a general approach to the
symbolic computation of the asymptotic bias and variance to order
$O(1/n)$ for those estimators.  The definition of (1) or (2)
allows an index set $\{m,k,i \}$ to represent an estimator so that
Theorem 1 or 2 can be easily implemented symbolically.
\bigskip
\begin{center}
3. BIAS OF BURG ESTIMATORS IN AR(2)
\end{center}
\medskip
The stationary second-order autoregressive model may be written as
$z_t = \phi_1 z_{t-1} + \phi_2 z_{t-2} + a_t$, where $a_t$ are
normal and independently distributed with mean zero and variance
$\sigma^2$ and parameters $\phi_1$ and $\phi_2$ are in the
admissible region, $|\phi_2| < 1$, $\phi_1 + \phi_2 <1$ and
$\phi_2-\phi_1<1$. The Burg estimate for $\phi_2$ may be obtained
directly from Percival and Walden (1993, eqn. 416d) and then the
estimate for $\phi_1$ may be obtained using the Durbin-Levinson
algorithm. After simplification, these estimates may be written as
\begin{equation}
\hat \phi_2 = 1 - {C D^2 - 2 E D^2 \over C D^2 + 8 F^2 G - 4 F H D},
\hat \phi_1 = {2 F  \over D} (1-\hat \phi_2)
\newcounter{BurgARTwoB}
\setcounter{BurgARTwoB}{\value{equation}}
\end{equation}
where
\begin{eqnarray*}
C = {1 \over n-2} \sum\limits_{t=3}^n (z_t^2 + z_{t-2}^2),
D =  {1\over n-1}\sum\limits_{t=2}^n (z_t^2 + z_{t-1}^2),
E = {1\over n-2} \sum\limits_{t=3}^n (z_t^2 z_{t-2}^2),\\
F =  {1\over n-1}\sum\limits_{t=2}^n (z_t z_{t-1}),
G = {1\over n-2} \sum\limits_{t=3}^n z_{t-1}^2,
H = {1\over n-2} \sum\limits_{t=3}^n (z_t z_{t-1}+z_{t-2}z_{t-1}).
\end{eqnarray*}

Using a Taylor series expansion of $\hat \phi_1$ and $\hat \phi_2$
about $\mu_{\cal A} = \E(\cal A)$, where ${\cal A}=C,D,E,F,G\ {\rm
and}\ H$, the bias to order $O(1/n)$ may be expressed in terms of
asymptotic expectations of products and cross products involving
$C,D,E,F,G$ and $H$. There are six squared terms and fifteen cross
product terms involved in each expansion, that is, it is required
to compute and simplify for each of these twenty one expansions
involving $C,D,E,F,G$ and $H$. These terms may all be written in
terms of the unbiased estimate of the autocovariance, ${\cal
S}_{m,k,i}$.
%\begin{equation}
%{\cal S}_{m,k,i} = {1\over n+1-i} \sum\limits_{t=i}^n z_{t-m} z_{t-k}.
%\newcounter{Smki}
%\setcounter{Smki}{\value{equation}}
%\end{equation}
The required asymptotic expectations of each term in the
expansions are obtained by Theorem 1, that is,
\begin{equation}
\lim_{n \to \infty} n \Cov( {\cal S}_{m,k,i}, {\cal S}_{f,g,j} ) =
\sum_{ h=-\infty }^{\infty} T_h,
\newcounter{SmkiCovariance}
\setcounter{SmkiCovariance}{\value{equation}}
\end{equation}
where $T_h = \gamma(h) \gamma(h-p+q)+\gamma(h+q) \gamma(h-p)$,
$p=m-k$, $q=f-g$ and
\begin{equation}
\gamma(h) =  \frac{{{{\zeta }_2}}^{1 + h} - {{{\zeta }_1}}^2\,{{{\zeta }_2}}^{1 + h} +
    {{{\zeta }_1}}^{1 + h}\,\left(  {{{\zeta }_2}}^2 -1 \right) }{
    \left(  {{{\zeta }_1}}^2 -1 \right) \,
    \left( {{\zeta }_1} - {{\zeta }_2} \right) \,
    \left(  {{\zeta }_1}\,{{\zeta }_2} -1 \right) \,
    \left( {{{\zeta }_2}}^2 -1 \right) },
\newcounter{ACVF}
\setcounter{ACVF}{\value{equation}}
\end{equation}
where $h \ge 0$, $\zeta_1$ and $\zeta_2$ are the roots, assumed
distinct, of the polynomial $\zeta^2 -\phi_1 \zeta -\phi_2= 0$.
The order $n^{-1}$ coefficient of the covariance expansion of
${\cal S}_{m,k,i}$ and ${\cal S}_{f,g,j}$ given in eqn.
({\theSmkiCovariance}) may be evaluated symbolically by defining
an operator of ${\cal S}_{m,k,i}$ and ${\cal S}_{f,g,j}$,
$\LCOV[\{m,k,i\}\{f,g,j\}]$. To illustrate this symbolic method
consider the evaluation of $\lim_{n \to \infty} n \Cov(2C, H)$
which is one of the twenty one order $n^{-1}$ expansion
coefficients involving $C,D,E,F,G$ and $H$ mentioned above. It may
be obtained by
\begin{eqnarray*}
\lim_{n \to \infty} n \Cov(2C,H)
&=&2\{\LCOV[(\{0,0,3\}+\{2,2,3\})(\{0,1,3\}+\{2,1,3\})]\}\\
&=&2\{\LCOV[\{0,0,3\}\{0,1,3\}]+\LCOV[\{0,0,3\}\{2,1,3\}]\\
& &+\LCOV[\{2,2,3\}\{0,1,3\}]+\LCOV[\{2,2,3\}\{2,1,3\}]\},
\end{eqnarray*}
since $C ={\cal S}_{0,0,3}+{\cal S}_{2,2,3}$, $H =
{\cal S}_{0,1,3}+{\cal S}_{2,1,3}$, and $\LCOV[\cdot]$ follows the
linearity and the distributive law.

After algebraic simplification the biases to order $O(1/n)$ are
found to be $\E(\hat \phi_1 - \phi_1) \doteq -(\zeta_1 +
\zeta_2)/n$ and $\E(\hat \phi_2 - \phi_2) \doteq (3 \zeta_1 \zeta_2
-1)/n$. More simply, in terms of the original parameters we have
\begin{equation}
\E(\hat \phi_1 - \phi_1) \doteq  -\phi_1/n
\newcounter{BiasOne}
\setcounter{BiasOne}{\value{equation}}
\end{equation}
and
\begin{equation}
\E(\hat \phi_2 - \phi_2) \doteq -(1 + 3 \phi_2)/n.
\newcounter{BiasTwo}
\setcounter{BiasTwo}{\value{equation}}
\end{equation}
We verified, using the same approach, that eqns. ({\theBiasOne})
and ({\theBiasTwo}) also hold for the case of equal roots of the
polynomial $\zeta^2 -\phi_1 \zeta -\phi_2= 0$.

For the stationary second-order autoregressive model with an
unknown mean, the Burg estimators can be written as the same ratio
function of the biased estimators of the autocovariances,
$\overline {{\cal S}}_{m,k,i}$, as was given in eqn.
({\theBurgARTwoB}). The symbolic approach is similar to the known
mean case, but includes one more inner product associated with the
biases of those autocovariance estimators, $\overline {{\cal
S}}_{m,k,i}$. The required asymptotic biases and covariances of
$\overline {{\cal S}}_{m,k,i}$ are obtained by Theorem 2.  The
biases to order $O(1/n)$ are found to be $\E(\hat\phi_1 - \phi_1)
\doteq ((\zeta_1 \zeta_2 -\zeta_1 -\zeta_2)-1)/n$ and $\E(\hat
\phi_2 - \phi_2) \doteq (4 \zeta_1 \zeta_2 -2)/n$. That is
\begin{equation}
\E(\hat \phi_1 - \phi_1) \doteq -(\phi_2+\phi_1+1)/n
%\newcounter{BiasOne}
%\setcounter{BiasOne}{\value{equation}}
\end{equation}
and
\begin{equation}
\E(\hat \phi_2 - \phi_2) \doteq -(2 + 4 \phi_2)/n.
%\newcounter{BiasTwo}
%\setcounter{BiasTwo}{\value{equation}}
\end{equation}

Once an estimator of a stationary time series is written as a well
defined function composed of ${\cal S}_{m,k,i}$ or $ \overline
{{\cal S}}_{m,k,i}$, by expanding it by a Taylor series, the
estimate bias and variance to order $n^{-1}$ may be obtained by
Theorem 1 or 2 with symbolic computation.  This approach can be
applied in the bias derivation of the Burg estimator,
$\hat{\rho}$, in the first order autoregressive model, AR(1).  We
have obtained that its bias to order $n^{-1}$ is $-2\rho/n$ in a
zero mean case and $-(1+3\rho)/n$ in an unknown mean case.
Therefore, for both of AR(1) and AR(2) cases, the biases to order
$n^{-1}$ of the Burg estimators are the same as the least squares
estimation for a known mean case as well as for an unknown mean
case.  These results are consistent with the simulation study
reported by Tj\o stheim and Paulsen (1983).

\bigskip
\begin{center}
4. CONCLUDING REMARKS
\end{center}
\medskip
In addition to deriving the bias for the Burg estimator, we used
our computer algebra method to verify the bias results reported by
Tj\o stheim \& Paulsen (Correction, 1984) and we also carried out
simulation experiments which provided an additional check on our
results (Zhang, 2002). {\it Mathematica\/} (Wolfram, 1999)
notebooks with the complete details of our derivations are
available on request from the authors.

Since many quadratic statistics in a stationary time series can be
expressed in terms of ${\cal S}_{m,k,i}$ or $ \overline {{\cal
S}}_{m,k,i}$, our computer algebra approach can be applied to
derive their laborious moment expansions to order $O(1/n)$.  As examples, using our
method, we can easily obtain the results by Bartlett (1946),
Kendall (1954), Marriott and Pope (1954), White (1961) and Tj\o
stheim and Paulsen (1983).

%\newpage
\bigskip
\begin{center}
{\bf REFERENCES\hfill}
\end{center}
\parindent 0pt

\hind ANDREWS, D. F. and STAFFORD, J. E. (2000) {\it Symbolic
Computation for Statistical Inference,\/} Oxford University Press.

\hind ANDREWS, D. F. and STAFFORD, J. E. (1993) Tools for the
Symbolic Computation of Asymptotic Expansions. {\it J. R. Statist.
Soc. B\/} {\bf 55, No. 3,\/} 613--627.

\hind BARTLETT, M. S. (1946)  On the Theoretical Specification and
Sampling Properties of Autocorrelated Time-Series. {\it J. R.
Statist. Soc. Suppl. \/} {\bf 8,\/} 27--41.

COOK, P. and BROEMELING, L.D. (1995)
Bayesian Statistics Using Mathematica. {\it The American Statistician\/}
{\bf 49\/} 70--76

\hind CURRIE, I.D. (1995)
Maximum Likelihood Estimation and Mathematica.
{\it Applied Statistics\/} {\bf 44,\/} 379--394.

\hind FULLER, W. A. (1996) {\it Introduction to Statistical Time
Series\/}, 2nd Ed., New York: Wiley.

\hind KENDALL, M. G. (1954) Notes on the bias in the estimation of
autocorrelation. {\it Biometrika \/} {\bf 41,\/} 403--404.

\hind
MARRIOTT, E. H. C. and POPE, J. A. (1954) Bias in the
estimation of autocorrelation. {\it Biometrika\/} {\bf 41,\/}
390-402.

\hind
LYSNE, D. and TJ\O STHEIM, D. (1987)
Loss of Spectral Peaks in Autoregressive Spectral Estimation.
{\it Biometrika\/} {\bf 74,\/} 200--206.

\hind
PERCIVAL, D. B. and WALDEN A. T. (1993)
{\it Spectral Analysis For Physical Applications,\/}
Cambridge: Cambridge University Press.

\hind
ROSE, C. and SMITH, M.D. (2002)
{\it Mathematical Statistics with Mathematica\/},
New York: Springer-Verlag.

\hind
ROSE, C. and SMITH, M.D. (2000)
Symbolic maximum likelihood estimation with Mathematica.
{\it The Statistician\/}, {\bf 49\/}, 229--240.

\hind SMITH, B. and FIELD, C. (2001) Symbolic Cumulant
Calculations for Frequency Domain Time Series.
{\it Statistics and Computing\/} {\bf 11}  75--82.

\hind TJ\O STHEIM, D. and PAULSEN, J. (1983) Bias of some
commonly-used time series estimates. {\it Biometrika\/} {\bf
70,\/} 389--399. Correction {\it Biometrika\/} {\bf 71,\/} p. 656.

\hind
WOLFRAM, S. (1999)
{\it The Mathematica Book\/}, 4th Ed.
Cambridge: Cambridge University Press.

\hind WHITE, J. S. (1961) Asymptotic Expansions for the Mean and
Variance of the Serial Correlation Coefficient.
{\it Biometrika\/} {\bf 48, \/}, p. 85--94.

\hind
ZHANG, Y. (2002)
{\it Topics in Autoregression\/},
Ph.D. Thesis, University of Western Ontario.

\end{document}